\documentclass[aps,showpacs,twocolumn,groupedaddress]{revtex4-2}

\usepackage{mathtools}
\usepackage{amsmath}
\usepackage{amssymb}
\usepackage{graphicx,epstopdf,subfigure}
\usepackage{xcolor}
\usepackage[normalem]{ulem}
\usepackage[capitalise]{cleveref}
\usepackage{tikz,bm,standalone}

\newcommand{\ie}{\emph{i.e.}}
\newcommand{\eg}{\emph{e.g.}}

\usepackage{color}

\newcommand{\as}[1]{[\textcolor{green!50!black} {AS: {#1}}]}

\begin{document}

\title{Network and Phase Symmetries\\
Reveal That\\
Amplitude Dynamics Stabilize Decoupled Oscillator Clusters}

\author{Jeffrey Emenheiser}
\email{jemenheiser@ucdavis.edu}
 
\author{Anastasiya Salova}
\email{avsalova@ucdavis.edu}

\author{Jordan Snyder}
\email{jsnyd@uw.edu}

\author{James P. Crutchfield}
\email{chaos@ucdavis.edu}

\author{Raissa M. D'Souza}
\email{raissa@cse.ucdavis.edu}
\affiliation{Complexity Sciences Center and\\
Physics, Mathematics, and Computer Science Departments,\\
University of California, Davis, One Shields Avenue,
Davis, California 95616}
\date{\today}

\begin{abstract}
Oscillator networks display intricate synchronization patterns. Determining
their stability typically requires incorporating the symmetries of the network coupling. Going beyond analyses that appeal only to a network's
automorphism group, we explore synchronization patterns that emerge from the
phase-shift invariance of the dynamical equations and symmetries in the nodes.
We show that these nonstructural symmetries simplify stability calculations.
We analyze a ring-network of phase-amplitude oscillators that exhibits a
``decoupled" state in which physically-coupled nodes appear to act
independently due to emergent cancellations in the equations of dynamical
evolution. We establish that this state can be linearly stable for a ring of
phase-amplitude oscillators, but not for a ring of phase-only oscillators that
otherwise require explicit long-range, nonpairwise, or nonphase coupling. In
short, amplitude-phase interactions are key to stable synchronization at a
distance.
\end{abstract}

\maketitle

\section{Overview} 

Oscillator networks exhibit a wide variety of coordinated and collective
behaviors, \eg, globally synchronized states, splay states, cluster states, and
chimera states \cite{belykh2008cluster, cho2017stable,
nicolaou2019multifaceted, kuramoto1975, brown2003globally,
strogatz2000kuramoto, abrams2004chimera, ashwin1992dynamics,
pecora2014cluster}. Stability calculations for cluster synchronization
\cite{pecora2014cluster,sorrentino2016complete} and
independently-synchronizable clusters \cite{cho2017stable} were recently
simplified using the underlying structural symmetries of the network connecting
the oscillators. Yet, in addition to network connectivity, general symmetries
may play a significant role in collective behavior. It is well known, for
instance, that general dynamical systems exhibit emergent symmetries
\cite{golubitsky2003symmetry}. We adapt this
insight to show that dynamical symmetries---phase-shift symmetries in nodal
dynamics and coupling---are key to determining oscillator network stability.
Most notably, these symmetries reveal collective states with nuanced
behaviors---collective states that \emph{synchronize at a distance}---whose
stability properties can be understood using dynamical symmetries.

We focus, in particular, on a set of trajectories dubbed the ``decoupled" state
in a ring-network of phase-amplitude oscillators---a collective state whose
existence was conjectured some time ago \cite{alexander1988global, emenheiser2016patterns}, but only recently realized experimentally
\cite{matheny2019exotic}. This synchronization at a distance arises when
oscillators, physically coupled along a nearest-neighbor ring network, begin to
act independently of their immediate neighbors due to emergent cancellations in
the dynamical equations of motion. The result is collective dynamical patterns
that are higher-order and longer-range than the pairwise oscillator
physical coupling. For instance, the decoupled state experimentally explored in
Ref.~\cite{matheny2019exotic} exhibits effective next-nearest-neighbor coupling
in a network formed from purely nearest-neighbor physical coupling.

The decoupled state generally appears through symmetry breaking. To analyze
stability we use the fact that, for each symmetry-breaking state, there is a
corresponding subgroup of the system's full symmetry group that leaves each
component of that state invariant. The Jacobian, used to monitor the state's
stability, must commute with the subgroup symmetry operators at each point in
the state's trajectory \cite{golubitsky2003symmetry}. This, in effect, defines
``symmetry-breaking'' by considering those symmetries of the dynamics that are
not broken. Mutually diagonalizing these subgroup symmetry operators determines
the corresponding block diagonal form of the Jacobian. If the block size is
sufficiently small, we give the Jacobian's eigenvalues and eigenvectors in
closed-form. This tool allows us to investigate the stability of the decoupled
state, even with arbitrarily large ring networks.

In this way, we establish that the decoupled state can be linearly stable for a
ring of phase-amplitude oscillators, but not for a ring of phase-only
oscillators with phase-based coupling (\eg, as in the Kuramoto model
\cite{kuramoto1975,kuramoto1984book}). This demonstrates the importance of
oscillator amplitude degrees of freedom for stability and for creating
long-range effective couplings.

Our development proceeds as follows. First, we describe the oscillator network
and its symmetry subgroups and show how these constrain the Jacobian. Second,
we introduce a specific dynamic corresponding to the system experimentally
studied in Ref. \cite{matheny2019exotic} consisting of a ring of
nanoelectromechanical, phase-amplitude oscillators. Third, we show that
symmetry considerations predict the existence of the decoupled state (similar
to Refs. \cite{alexander1988global, ashwin1992dynamics}) and, moreover,
simplify the stability analysis. Fourth, when the decoupled state was achieved
in experiment, drift in the phase difference between the decoupled clusters was
observed \cite{matheny2019exotic}. We explain this by showing that small
deviations in the oscillators' natural frequencies break the symmetry and
result in drift. Finally, with drift, the Jacobian becomes time dependent and
so we use Floquet theory to analyze network stability. This analysis shows that
alternating the natural frequencies of adjacent oscillators introduces an intricate
dependence of the stability on natural frequency difference and other system
parameters.

\section{Symmetries and Stability}
\label{sec:symmetry}

Consider a generic continuous-time dynamical system $\dot{x}=f(x)$ defined on a
state space $X$. A set of invertible operations $\{\gamma: X\to X\}$, together
with operator composition, generates a symmetry group $\Gamma$ of the dynamics
if and only if, for all $x\in X$, each operation obeys:
\begin{align}
	\label{eq:symmetry}
f(\gamma x) = (D_x\gamma)f(x)
  ~,
\end{align}
where $D_x\gamma\colon T_x X \to T_x X$ is the differential operation at $x$.
This defines an \emph{equivariant} dynamical system with respect to the
symmetry group \cite{moehlis2007equivariant}. When the tangent ($T$) or
differential ($D$) operators are independent of reference point $x$, we drop
the subscript.

Under the restriction that each $\gamma$ acts linearly on $X$, the differential
operation becomes trivial, and the equivariance condition becomes: $f(\gamma x)
= \gamma f(x)$. This means one can side-step differences between $X$ and $TX$.
While such conditions are typically met by symmetry operations associated with the network structure, the
nonlinear operations here require more careful handling. Similarly, curvilinear
coordinates may disrupt the matrix representations of linear operations. In
these cases, one must be mindful of the symmetry differences between the state
space $X$ and tangent space $T X$.

Given a dynamical system that respects a symmetry group $\Gamma$, the
\emph{isotropy subgroups} $\Sigma \le \Gamma$ are those that leave state-space
subsets invariant. A particular trajectory in state space, such as a fixed
point or limit cycle, is said to respect an isotropy subgroup $\Sigma$ if all
points $x(t)$ in the trajectory are invariant under $\Sigma$:
\begin{align*}
\sigma x(t) = x(t)
  ~,
\end{align*}
for all $\sigma \in \Sigma$ and all $t \in \mathbf{R}$. For a given dynamical
system, classes of trajectories may be predicted by identifying the system's
largest symmetry group and then identifying subgroups $\Sigma$ that leave a
nontrivial set of points in state space fixed. These state-space subsets are
themselves dynamically invariant. 

$\Sigma$'s group structure provides a convenient coordinate basis for
describing the evolution of states close to the trajectory of interest and,
therefore, for understanding its stability. In particular, we study the
evolution of infinitesimal deviations $x(t) + \epsilon\delta x(t)$, $\epsilon
\ll 1$. To leading order in $\epsilon$, their evolution is governed by the
linear ordinary differential equation: 
\begin{align*}
\frac{d}{dt} \delta x = J(x) \delta x
  ~,
\end{align*}
where $J = \frac{\partial f}{\partial x}$ is the system Jacobian.

Since the global evolution respects the symmetry operators $\sigma$, the linear
dynamics respects symmetry operators $D\sigma(x)$. This means that the Jacobian
and the symmetry operators commute at each point in the fixed-point subspace of
$\Sigma$:
\begin{align*}
D\sigma(x)J(x) = J(x)D\sigma(x)
  ~.
\end{align*}
The Jacobian thus shares eigenspaces with each of the differential group
operators. To find these eigenspaces, we block diagonalize the matrices
corresponding to the symmetry group generators. This can be done by finding the
isotypic components of the differential group $\{D \sigma: \sigma \in
\Sigma\}$: each block corresponds to an irreducible representation of the group
\cite{golubitsky2003symmetry}. Since they commute, the linear operators
$D\sigma(x)$ and $J(x)$ then share $\Sigma$-irreducible invariant subspaces
\cite{golubitsky2012singularities}, acting on each subspace according to the
corresponding diagonal block.


Below, we show how a behavior's symmetry may be used to block-diagonalize the
linear dynamics around any point of that trajectory class and so determine its
stability properties. Assessing this class then reduces to two distinct
problems: (i) evolving states within the subset of state space invariant to
a known symmetry subgroup and (ii) evolving perturbations transverse to that
subspace.


\section{Phase-amplitude Oscillator Ring}
\label{sec:dynamics}

We now turn to the specific system under study---a ring of
nanoelectromechanical (NEMS) phase-amplitude oscillators. For the experiments
reported in Ref. \cite{matheny2019exotic} each oscillator is implemented as a
piezoelectric oscillating membrane that is coupled to other oscillators
electronically in a tunable network. The resulting oscillator networks
exhibited a wide variety of exotic synchronization patterns including the
decoupled state, mentioned above.

In fact, as found in Ref. \cite{matheny2019exotic}, for a ring network the
dynamical evolution is captured by well-understood equations of motion
\cite{lifshitz2008nonlinear, villanueva2013surpassing}. Representing the state
of each of the $N$ oscillators as a complex number $A_j$, with node index $j =
0, 1, ..., N-1$, the dynamics are:
\begin{align}
  \frac{dA_j}{dt} = -A_j + i\omega_jA_j &+ 2i\alpha|A_j|^2A_j + \frac{A_j}{|A_j|}\nonumber\\
  &+ i\beta\left(A_{j-1}-2A_j + A_{j+1}\right)
  ~.
\label{eq:complexphase}
\end{align}
Indexing is taken \textit{modulo} $N$. For symmetry considerations discussed
later, we limit our development to rings consisting of multiples of four
($N=4M$) oscillators.

There are three categories of tunable system parameters: $\omega_j$ is the natural frequency
of oscillator $j$, $\alpha$ controls the Duffing nonlinearity (equal across
all oscillators), and $\beta$ is the reactive coupling strength between
adjacent oscillators (equal across all pairs of neighbors). The dynamics of an individual oscillator in
\cref{eq:complexphase} is similar to that of the widely-studied Stuart-Landau
oscillator.

It is convenient to represent the state in $\mathbf{R}^{2N}$, where each
oscillator's state is split into real amplitude and phase components:
$a_je^{i\phi_j} = A_j$. The equations of motion become:
\begin{align}
  \label{eq:amp}
  \frac{da_j}{dt} & = 1-a_j - \beta a_{j-1}\sin\left(\phi_{j-1}-\phi_j\right)
  \nonumber \\
  & \qquad - \beta a_{j+1}\sin\left(\phi_{j+1}-\phi_j\right) \\
  \label{eq:phase}
  \frac{d\phi_j}{dt} & = \omega_j + 2\alpha a_j^2
  + \beta\frac{a_{j-1}}{a_j}\cos\left(\phi_{j-1}-\phi_j\right)  \nonumber \\
  & \qquad + \beta\frac{a_{j+1}}{a_j}\cos\left(\phi_{j+1}-\phi_j\right)-2\beta 
  ~.
\end{align}

We now consider this system's symmetries; \ie, we identify the operations
that satisfy \cref{eq:symmetry}.

If all oscillators have \emph{uniform} frequency, $\omega_j=\omega$, the ring
dynamics respects the symmetry group generated by rotations $\sigma_\text{rot}:
A_j \mapsto A_{j+1}$ and by (node-centered) reflections $\sigma_\text{ref}: A_j
\mapsto A_{N-j}$ of the ring. The symmetry of the undirected cycle is called
the \emph{dihedral group} and denoted $\mathcal{D}_N$. Its elements are the
$2N$ unique products of $\sigma_\text{rot}$ and $\sigma_\text{ref}$. Elements
of the dihedral group act as permutation matrices on both the complex
coordinates of \cref{eq:complexphase} and the real coordinates of
\cref{eq:amp,eq:phase}, as is standard for topological symmetries. These
operations merely reorder the nodal coordinates.

We also consider the case of \emph{alternating} frequencies: $\omega_j =
\omega\mp\Omega/2$, where even numbered $j$ follow the $-$ and odd numbered $j$
the $+$, and we call $\Omega$ the \emph{detuning}. (We will see that the mean
frequency $\omega$ does not influence stability calculations and merely defines
a frame of reference.) The network rotational symmetry is now reduced to that
generated by $\sigma_{\text{rot}}^2$. The reflectional symmetry
$\sigma_\text{ref}$ is still respected, generating a $\mathcal{D}_{N/2}$
symmetry group.

Regardless of node frequency details, the ring also respects a continuous
symmetry of uniform phase shifts: $\sigma_\theta: \{\phi_j\} \mapsto \{\phi_j +
\theta\}$ with any $0 \le \theta < 2\pi$. We denote this continuous group
$\mathcal{T}$. In real amplitude and phase coordinates, this operation is
affine yet nonlinear. Its differential $D_x\sigma_\theta$ is equivalent to the
identity at all points $x$. This can be seen in \cref{eq:phase,eq:amp} by
shifting all phases by $\theta$: these phase shifts cancel and the equations of
motion are invariant. However, in the complex amplitudes of
\cref{eq:complexphase}, this becomes the linear action $\sigma_\theta: \{A_j\}
\mapsto \{e^{i\theta} A_j\}$.

Since uniform phase shifts commute with reordering node indices, the full symmetry group of the system is the direct product group $\mathcal{D}_N \times \mathcal{T}$ in the case of uniform frequencies (and $\mathcal{D}_{N/2} \times \mathcal{T}$ in the case of alternating frequencies). Subgroups of $\mathcal{D}_N \times \mathcal{T}$ may contain nontrivial phase
shifts, revealing interesting synchronization patterns, as framed generally
by Ref. \cite{ashwin1992dynamics} and used in Ref. \cite{matheny2019exotic} to characterize synchronization patterns in the NEMS system.

For each subgroup of system symmetries, the set of points left unchanged by the
action of all symmetry operators defines an invariant set of the system
dynamics. Time-evolution cannot break a symmetry of the initial condition, if
the dynamics themselves respect that symmetry. To proceed, we identify a
particular symmetry subgroup that defines an interesting invariant set in both
the uniform and alternating frequency cases. 



\begin{figure*}[t]
\includegraphics{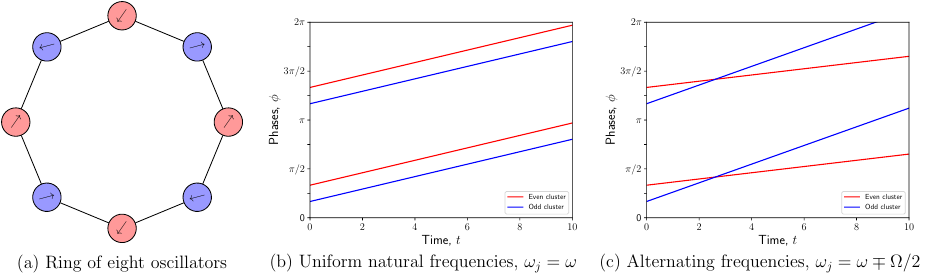}
\caption{(a) Snapshot of a ring of eight oscillators in the decoupled state,
	with nodes colored according to synchronized cluster and arrows indicating
	local phase. Solid black lines indicate physical coupling between 
	oscillators.
	(b)~Reference phase of each decoupled cluster versus time for oscillators
	with uniform natural frequencies,
	demonstrating constant phase differences.
	(c)~Reference phase of each decoupled cluster versus time for oscillators
	with alternating natural frequencies, demonstrating drift in phase
	difference between nodes in different clusters.
	The nonlinearity, coupling strength, and mean natural frequency are
	$\alpha=0.1, \beta=1.0,$ and $\omega=2$, respectively. The difference in
	natural frequencies between the two clusters is
	(b):~$\Omega=0$ and (c):~$\Omega=0.2$.
	}
\label{fig:behavior}
\end{figure*}

\section{Decoupled Antiphase-Synchronized Clusters}
\label{sec:behavior}

With this background in mind, we turn to the collective state of interest---the
decoupled state in which each node appears to act independently of neighbors to
which it is physically coupled. For the equations of motion in
\cref{eq:complexphase}, independence is the condition that $A_{j-1} + A_{j+1} =
0$; that is, next-nearest neighbors must be locked in perfect antiphase. This
condition causes the coupling terms to cancel out, so that each oscillator
appears to evolve under the influence of itself alone.  On a ring of $N$
oscillators, this is only satisfied if $N=4M$. The system then splits into
a group of antiphase synchronized even-numbered oscillators and a group of
antiphase synchronized odd-numbered oscillators, as illustrated in
Fig.~\ref{fig:behavior}(a). Remarkably, there is no constraint on the phase
differences between the two groups.

This decoupled state reflects a system symmetry. Note that the requirement
$A_{j-1} + A_{j+1} = 0$ is exactly the condition for a state to be invariant to
the operation $\sigma_\pi\sigma_\text{rot}^2$. This operator rotates the ring
by two oscillators and advances the phase of all oscillators by one half
period: $\sigma_\pi\sigma_\text{rot}^2A_j = -A_{j+2}$. It generates a cyclic
group of order $N/2$ and is a member of the symmetries of both the uniform
frequency and alternating frequency cases. The decoupled states, such as the
one shown in Fig.~\ref{fig:behavior}(a), lie precisely in the fixed-point
subspace of this subgroup of the system symmetry and, therefore, are an
invariant set of the dynamics.

This requirement also simplifies the coupling term in Eq.
(\ref{eq:complexphase}) ({\it i.e.}, the final term) to $-2i\beta A_j$. This
means that the evolution of each oscillator is determined only by its own
dynamic state. Using this in Eqs.  (\ref{eq:amp}) and (\ref{eq:phase}) one sees
that each oscillator amplitude will approach unity, behaving as if it were
uncoupled. 

An oscillator's exact phase is arbitrary, as is the phase difference between
the even- and odd-numbered oscillators. Defining a global reference phase
$\theta$ and phase difference $\psi$, we obtain the set of solutions describing
the decoupled state:
\begin{align}
a_j(t) &= 1, \nonumber \\
\phi_j(t) &= \left\{\begin{array}{lr}
	\theta(t)&\\
	\theta(t) + \psi(t)&\\
	\theta(t) + \pi&\\
	\theta(t) + \psi(t) + \pi,&
	\end{array}
	\quad j\!\!\mod 4 = \begin{array}{c}0\\1\\2\\3\end{array}\right. 
  ~.
\label{eq:torus}
\end{align}
For both the uniform- and  alternating-frequency cases:
\begin{align*}
d\theta/dt & = \omega + 2\alpha - 2\beta  \\
  d\psi/dt & = \Omega
  ~.
\end{align*}

Recall that $\Omega=0$ for uniform frequencies and $\Omega \ne 0$ for
alternating frequencies, so $\psi$ is fixed in time in the former case, but
drifts according to the detuning $\Omega$ for the latter. Examples are shown in
Fig.~\ref{fig:behavior} (b) and (c). The initial values of $\theta$ and $\psi$
are free variables defining a $2$-torus. The solutions over all possible
initial values form a set of limit cycles that foliate the torus.


\section{Stability}
\label{sec:static}

This particular decoupled state (with multiple anti-synchronized clusters) was
previously studied via symmetry considerations, for instance, on weakly coupled
identical phase oscillators \cite{ashwin1992dynamics, brown2003globally} and
linearly coupled oscillators \cite{alexander1986global,alexander1988global}.
Stability properties were not addressed. Recently, it was shown how a
broad variety of decoupled states can arise from balanced equivalence relations
rather than from symmetries and that symmetries are only required in an
effective network of clusters of nodes \cite{salova2020decoupled}. However,
Ref. \cite{salova2020decoupled} addressed only the stability properties of
networks of uniform nodes. Here, we provide stability analysis in cases of both
uniform and alternating frequencies and for a specific dynamics.

\begin{figure*}[t]
\includegraphics{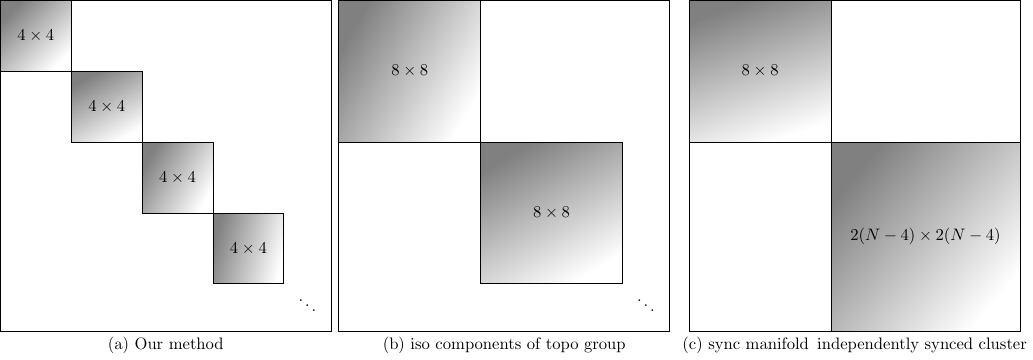}
\caption{Guaranteed zero block structure of Jacobian matrix in various coordinates. Non-zero entries may occur only in the grayscale blocks. The coordinates are given by the following methods:
	(a) Isotypic components of $\sigma_{\text{rot}}^2\sigma_\pi$ (via method introduced here).
	(b) Isotypic components of $\sigma_{\text{rot}}^4$ (via method
	of Ref. \cite{sorrentino2016complete}), and
	(c) Synchronization manifold and independently synchronized cluster set
	(via method of Ref. \cite{cho2017stable}).
	}
\label{fig:blocks}
\end{figure*}

To set up the linear stability analysis of the decoupled state, we next use the
symmetry structure to block-diagonalize the Jacobian. Given the solution in \cref{eq:torus} to the dynamics of
\cref{eq:amp,eq:phase}, the linear dynamics of
real deviations to amplitudes $\delta a$ and phases $\delta \phi$ are:
\begin{align}
\frac{d\delta a_j}{dt} & = \bigg[\beta\sin\psi\delta a_{j-1}
  - \delta a_j - \beta\sin\psi\delta a_{j+1} \nonumber\\
  & \qquad \mp \beta\cos\psi\delta \phi_{j-1}
  \pm \beta\cos\psi\delta \phi_{j+1}\bigg] \\
\frac{d\delta \phi_j}{dt} & = \bigg[ \pm\beta\cos\psi\delta a_{j-1}
  + 4\alpha\delta a_j \mp\beta\cos\psi\delta a_{j+1} \nonumber\\
  & \qquad + \beta\sin\psi\delta \phi_{j-1}
  -\beta\sin\psi\delta \phi_{j+1}\bigg]
  ~,
\label{eq:linearization}
\end{align}
where perturbations to even- (odd-)numbered oscillators follow the upper
(lower) $\pm/\mp$ option. This system of linear ODEs is independent of the
global phase $\theta$ and depends on the phase difference $\psi$.  
 
We simplify this linearization by finding convenient coordinates. With the
solution contained in the fixed-point subspace of the cyclic group $\Sigma$ generated by $\sigma_\text{rot}^2\sigma_\pi$, we know that $D\sigma J = J D\sigma$ for all $\sigma\in
\Sigma$. Equivalently through $\{D\sigma\}$'s isotypic components and
through eigendecomposition of the generating $D\sigma$, we find a coordinate
system in which the Jacobian is block diagonal, consisting of $N/2$ blocks of
size $4\times 4$. First, we define $\zeta^k = e^{(4k i \pi/N)}$ to be the
(N/2)th root of unity. This corresponds to the $N/2$ order of our cyclic group.
Let us define vectors of length $N/2$: $\Phi^{(k)}_j = \zeta^{k\cdot j}$, with
$k=0, 1, ..., \tfrac{N}{2}-1$. These $\Phi^{(k)}$ vectors are wave patterns
over $N/2$ elements. Each element of a $\Phi^{(k)}$ is associated with one
adjacent pair of oscillators, each with two real degrees of freedom. We
therefore take the matrix outer product with two $2\times 2$ identity matrices,
$I_2$: the first for the pair of oscillators and the second for the amplitude
and phase of each. This creates the $2N \times 4$ matrices $V^{(k)} =
I_2\otimes\Phi^{(k)}\otimes I_2$. These matrices, in fact, define $4$ eigenvectors of the generating (and therefore every) group operator. The columns share an eigenvalue: $\sigma_\text{rot}^2\sigma_\pi V^{(k)} = \zeta^k V^{(k)}$.

Since $V^{(k)}$'s columns exactly span an eigenspace (isotypic component), the Jacobian may be written in a coordinate system in which it is block diagonal with $2M$ blocks of size $4\times 4$, with each block specified by $D_k = V^{(k)T} J V^{(k)}$:
\begin{widetext}
\begin{align}
  \label{eq:blocks}
  &D_k = \frac{1}{2}
  \left[\begin{array}{cccc}
    -1 & -\beta(1-\zeta^{-k})\sin\psi & 0 & \beta(1-\zeta^{-k})\cos\psi \\
    \beta(1-\zeta^k)\sin\psi & -1 & \beta(1-\zeta^{k})\cos\psi & 0 \\
    4\alpha & -\beta(1-\zeta^{-k})\cos\psi & 0 & -\beta(1-\zeta^{-k})\sin\psi \\
    -\beta(1-\zeta^k)\cos\psi & 4\alpha & \beta(1-\zeta^k)\sin\psi & 0
  \end{array}\right]
  ~.
\end{align}
\end{widetext}
Importantly, this coordinate system is independent of system state, within the
solution set of interest (\cref{eq:torus}). The $V$'s used to compute each
block consist of constants, and the resulting Jacobian block depends only on
the phase difference $\psi$. Figure \ref{fig:blocks} shows the block structure,
along with the coarser block structure predicted by alternate methods that
consider purely-structural symmetries.

Instead of phase-amplitude oscillators, consider a ring of phase-only
oscillators. Then, the decoupled state specified in Eq. (\ref{eq:torus}) is a
guaranteed solution for $D_{4M}$ networks of identical nodes. Each has a $\mathcal{T}$
phase symmetry, so long as the coupling function respects a parity condition:
$g(\psi) + g(\pi-\psi) = 0$.  However, the diagonal elements of the Jacobian
are proportional to this term. This means that the Jacobian has zero trace.
Since the trace is the sum of eigenvalues, any linearly stable mode implies at
least one linearly unstable mode. And, so the decoupled state is not stable for phase-only oscillators.

Achieving a linearly-stable decoupled state for phase-only oscillators requires
introducing long-range, nonpairwise, or nonphase coupling. Thus, although the
amplitudes of the phase-amplitude oscillators have a fixed point at unity, the
amplitude degree of freedom plays a central role in stability.

\subsection{Uniform Frequencies: $\omega_j=\omega$}

With uniform frequencies, the phase difference $\psi$ between decoupled
clusters is constant and the Jacobian is fixed in time. The decoupled state's
stability is thus given by the real parts of the Jacobian's $2N$ eigenvalues.
This problem then reduces to finding eigenvalues of its $2M$ blocks each of
size $4\times 4$. This requires finding the roots of quartic polynomials, which
are:
\begin{widetext}
\begin{align}
\lambda = -\frac{1}{4} \pm \frac{1}{4}\sqrt{1 - 8\beta^2\left(1-\cos\frac{k\pi}{M}\right) \pm 4\beta \sqrt{2\left(16\alpha^2\cos^2\psi - \sin^2\psi\right)\left(1-\cos\frac{k\pi}{M}\right)}}
\label{eq:eigenvalues}
\end{align}
\end{widetext}
The eigenvalues of Eq. (\ref{eq:linearization})'s linear dynamics are given
by Eq. (\ref{eq:eigenvalues}), for $k = 0, 1, ..., 2M-1$ and with the two $\pm$
options taken independently.

Note that there are four possible, adjacent phase differences in Eq.
(\ref{eq:torus}): $\psi, \pi-\psi, \pi+\psi, 2\pi-\psi$. Eigenvalue $\lambda$'s
dependence on $\psi$ arises only through $\sin^2$ and $\cos^2$, which
necessarily are equal for all four possible phase differences. This supports
the physical equivalence of these values of $\psi$.

\begin{figure}[t]
\centering
\includegraphics[width=\columnwidth]{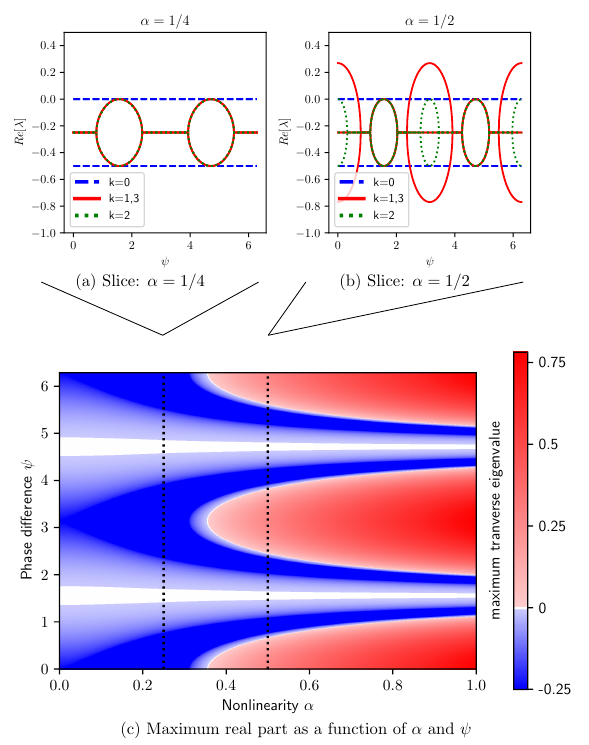}
\caption{Stability for uniform frequencies, $\omega_j = \omega$, at
	fixed coupling $\beta = 1$. Panels (a) and (b) show the real part of
	$D_k$'s eigenvalues, Eq.~(\ref{eq:eigenvalues}), versus phase difference
	$\psi$ for two choices of Duffing coefficient $\alpha$. Panel (c) shows the
	maximum real part of eigenvalues of blocks $D_{k\neq 0}$, quantifying the
	least stable perturbation that breaks the state's symmetry as a function of
	$\psi$ and $\alpha$.
	}
\label{fig:uniform}
\end{figure}

Figures~\ref{fig:uniform} (a) and (b) show the real part of $D_k$'s
eigenvalues, obtained from Eq. (\ref{eq:eigenvalues}) for $N=8$ and $\beta=1$,
as a function of phase difference $\psi$ and for different values of
nonlinearity: (a) $\alpha = 1/4$ and (b) $\alpha=1/2$. Each $k$-block gives
four eigenvalues, symmetric around $-1/4$ as expected from Eq.
(\ref{eq:eigenvalues}). The blue dashed lines at $\text{Re}(\lambda)=0$ show
neutral stability within the $k=0$ block. These eigenvalues are two-fold
degenerate and define the torus of solutions ({\it i.e.}, the class of
trajectories) of interest. The transverse perturbations, $k \neq 0$, exhibit
neutral stability at $\psi=\{\pi/2,3\pi/2\}$ at both values of
nonlinearity. The instabilities in the larger nonlinearity are centered at
$\psi = \{0,\pi\}$, where there is a second pair of regions in which the outer
square root of \cref{eq:eigenvalues} exhibits a real part.

Excluding the $k=0$ block, \cref{fig:uniform}(c) shows the maximum real part of
$D_k$'s eigenvalues for $N=8$ and $\beta=1$ as a function of both phase
difference $\psi$ and nonlinearity $\alpha$. This highlights the state's
largest instability at that point, with red being unstable, blue
stable, and white linearly neutral. Here, we see the steady neutral stability at
$\psi=\{\pi/2,3\pi/2\}$, but we also see the instabilities growing from
$\text{Re}(\lambda)= -1/4$ at $\psi = \{0,\pi\}$, as $\alpha$ increases.  From
\cref{fig:uniform}(b), we know that these unstable bands are in the $k=1,3$
block.

Due to the block diagonal structure induced by this state's symmetries, we can
precisely identify instabilities on the decoupled state: \cref{eq:eigenvalues}
directly relates parameter values with the growth and decay of perturbations.
All of the results presented in \cref{fig:uniform} are calculated as
closed-form expressions, capturing the eigenvalues of the $16\times 16$
Jacobian matrix.

\subsection{Alternating Frequencies, $\omega_j = \omega \mp \Omega/2$}
\label{sec:traveling}

In the experimental study of nanoelectromechanical  oscillator networks  
\cite{matheny2019exotic} parameters $\alpha$ and $\beta$ can be precisely
controlled. While the $\omega_j$'s can be tuned quite close to one another, small
deviations exist that break the $D_{4M}$ symmetry. When stable, all
perturbations transverse to the fixed-point subspace are exponentially
restored, and small dispersion in natural frequencies introduces a linear drift
within the invariant subspace. These new drifting states remain within the
space swept out by the uniform frequency states---sweeping $\theta$ and $\psi$
in \cref{eq:torus}---but support a drifting phase difference $\psi$.

We capture this behavior theoretically by alternating the natural frequencies
along the ring: $\omega_j = \omega\mp \Omega/2$. The form of Eq.
(\ref{eq:torus}) remains a solution to \cref{eq:complexphase}, but neighboring
phase differences now have a well-defined drift that depends on the magnitude
of the natural frequency difference between neighboring nodes: $d\psi/dt =
\Omega$.  Although the full system symmetry has changed, the solutions of
interest respect the same symmetries---the group generated by
$\sigma_\pi\sigma_\text{rot}^2$---as in the uniform frequency case. This leads
to the same block diagonalized Jacobian: \cref{eq:blocks}.

The Jacobian is now time-periodic. And, this requires Floquet theory for
stability analysis, evolving each of a set of vectors for one whole period. By
choosing this set as a basis for tangent space, we build the monodromy
matrix---the linear map corresponding to evolution of perturbations through a
period. Matrix eigenvalues capture how the perturbations evolve. The
associated Floquet exponents are the natural logarithm of the magnitude of its
eigenvalues, with positive values implying unstable growth and negative values
showing stable decay. If the Jacobian were constant (\ie, periodic with any
stated period), this procedure returns the real parts of the Jacobian
eigenvalues themselves, as desired.

We performed the integration using the standard fourth-order Runge-Kutta scheme
with $1,000$ integration steps~\cite{suli2003introduction}. The Jacobian's
diagonal form allows each $4\times 4$ block to be integrated independently,
significantly reducing the computational overhead for arbitrarily large
oscillator rings.

\begin{figure}[t]
\centering
\includegraphics[width=\columnwidth]{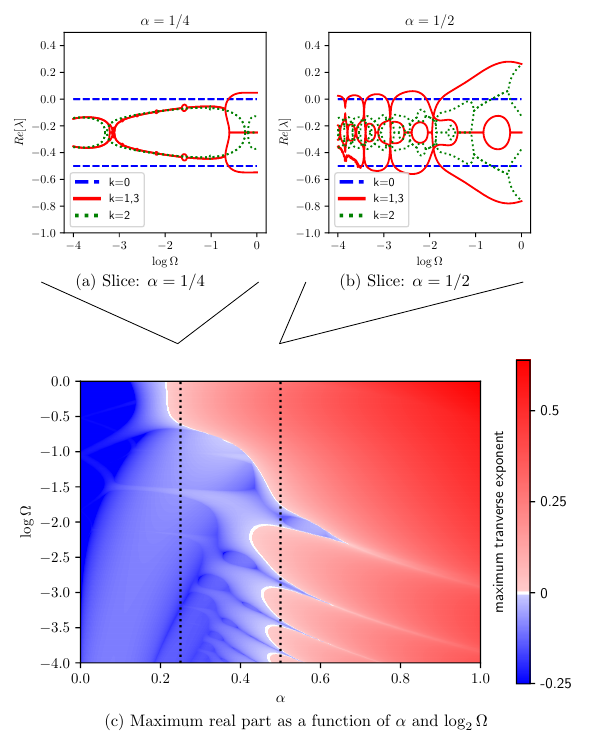}
\caption{Stability for alternating natural frequencies
	$\omega_j = \omega \mp \Omega/2$ with fixed coupling $\beta = 1$: Panels
	(a) and (b) plot the real part of the Floquet exponents of $D_k$
	(Eq.~\ref{eq:blocks}) versus relative frequency $\Omega$ for two choices of
	Duffing coefficient $\alpha$. Panel (c) shows the maximum real part of
	$D_{k\neq 0}$'s Floquet exponents, tracking the least stable perturbation
	that breaks the state's symmetry as a function of $\log_2 \Omega$.
	}
\label{fig:alternating}  
\end{figure}

We find that the stability of the resulting state depends intricately on system
parameters, including the natural frequency difference between adjacent
oscillators, as shown in Fig.~\ref{fig:alternating} \footnote{Noisy results
at low $\Omega$ in \cref{fig:alternating}(b) are a numerical artifact due to
long integration times and nonzero Floquet exponents}. Rather than being a
function of phase difference, stability is now a function of rate $\Omega$ at
which the system drifts through phase differences.

Here, we see many overlapping bands of instability emerging from
$\text{Re}(\lambda)=-1/4$. A series of such bands, going unstable around
$\alpha = 0.45$, appear to again come from the $k=1, 3$ blocks and leave small
windows of stability---narrow ranges of $\Omega$ in which the state
is stable.

Even though alternating frequencies in the oscillator ring break the system
symmetry and induce a time-periodic Jacobian, the symmetries of the state
itself are unchanged and the group-theoretic block-diagonalization is remains
valid. This allows the stability analysis to be performed on four dimensional
subspaces, even for arbitrarily large rings.

\section{Discussion}
\label{sec:discussion}

We demonstrated how phase-shift symmetries in the nodes and the dynamics, together with symmetries in
the network connectivity, can simplify stability calculations. This extends
recent results on cluster synchronization to include symmetries beyond the
node-connectivity automorphism group.  
We then used this approach to analyze the stability properties of attractors that emerge in a ring of phase-amplitude oscillators. We identified a dynamically decoupled
state as particularly novel, showing that it is stable for a ring a
phase-amplitude oscillators, but not linearly stable for the analogous ring of
phase-only oscillators.

Our results here highlight the
importance of oscillator amplitudes in generating polyadic and long-range
effective coupling. 
A reduced phase equation for the NEMS dynamics is introduced in~\cite{matheny2019exotic} 
but, to first order in coupling, it only includes the
nearest-neighbor phase interactions. This illustrates the need for future
analysis that either formally moves to higher-order phase models or retains
anharmonicity via amplitudes. 

Our methodology and results apply directly to a variety of similar models,
including any anharmonic oscillator with an attracting oscillation-amplitude,
such as Stuart-Landau oscillators \cite{salova2020decoupled}. The equivalence
of the decoupled state to a fixed-point subspace is, in such cases, a result
of linear coupling. That said, the coupling could have real rather than (or in
addition to) imaginary coefficient $i\beta \to K+i\beta$. This would directly
change the steady-state amplitudes rather than their frequencies. When $K>1$,
this leads to amplitude death. While the Jacobian block diagonalization is
ambivalent to coupling type, the precise stability results do not immediately
translate to systems where the coupling coefficient has a real component.

\begin{acknowledgments}

The authors thank Warren Fon, Mathew Matheny, Michael Roukes, Mehran Mesbahi, Afshin Mesbahi, and Leonardo Due{\~n}as-Osorio for helpful
discussions. This material is based upon work supported by, or in part by, the
U.S. Army Research Laboratory and the U. S. Army Research Office under MURI
award W911NF-13-1-0340 and grant W911NF-18-1-0028 and Intel Corporation support
of CSC as an Intel Parallel Computing Center.




\end{acknowledgments}

\bibliographystyle{unsrt}

\bibliography{gamma,chaos, biblio}

\end{document}